\documentclass[a4paper,10pt]{amsart}

\usepackage{amstext}
\usepackage{amsmath}
\usepackage{ifthen}
\usepackage{amscd}
\usepackage{amssymb}
\usepackage[french]{babel}
\usepackage[T1]{fontenc}
\usepackage[utf8]{inputenc}
\usepackage[all]{xy}
\usepackage{graphicx}
\usepackage{enumerate}
\usepackage{xspace}
\usepackage{pdfsync}
\usepackage{epic}
\usepackage{dsfont}

\newenvironment{demo}[1][]{\ifthenelse{\equal{#1}{}}{\noindent\textbf{D\' emonstration :}\xspace}{\noindent\textbf{D\' emonstration #1 :}\xspace}}{$\square$\newline}
\newtheoremstyle{break}% name
  {}%      Space above, empty = `usual value'
  {}%      Space below
  {\itshape}% Body font
  {}%         Indent amount (empty = no indent, \parindent = para indent)
  {\bfseries}% Thm head font
  {.}%        Punctuation after thm head
  {\newline}% Space after thm head: \newline = linebreak
  {}%         Thm head spec
  
\newtheoremstyle{rq}
  {}
  {}
  {\slshape}
  {}
  {\bfseries}
  {.}
  {3pt}
  {}

\newtheoremstyle{exemple}
  {}
  {}
  {\upshape}
  {}
  {\bfseries}
  {.}
  {3pt}
  {}

\newtheoremstyle{fact}
  {}
  {}
  {\slshape}
  {}
  {\bfseries}
  {.}
  {2pt}
  {}

\theoremstyle{fact}
\newtheorem*{fact*}{Affirmation}
\theoremstyle{break}
\newtheorem{thm*}{Th\' eor\`eme}
\newtheorem{thm}{Th\' eor\`eme}[section]
\newtheorem{conj}[thm]{Conjecture}
\newtheorem{cor}[thm]{Corollaire}
\newtheorem{lem}[thm]{Lemme}
\newtheorem{prop}[thm]{Proposition}
\newtheorem{defi}[thm]{D\' efinition}
\theoremstyle{rq}
\newtheorem{rem}[thm]{Remarque}
\theoremstyle{exemple}
\newtheorem{qt}[thm]{Question}
\newtheorem{ex}[thm]{Exemple}

\newcommand{\CC}{\mathbb C}
\newcommand{\QQ}{\mathbb Q}

\newcommand{\PP}{\mathbb P}

\newcommand{\kod}[1]{\kappa\left(#1\right)}

\newcommand{\To}{\longrightarrow}

\newcommand{\dimm}[1]{\mathrm{dim}(#1)}

\newcommand{\alb}{\operatorname{Alb}}

\newcommand{\rr}{\mathrm{r}}
\newcommand{\cc}{\mathrm{c}}
\newcommand{\diam}{\operatorname{diam}}
\newcommand{\longueur}{\operatorname{longueur}}

\title[Stabilit\' e des vari\' et\' es sp\' eciales]{Quelques propri\' et\' es de stabilit\' e des vari\' et\' es sp\' eciales}
\date{\today}
\author{Fr\' ederic Campana et Beno\^it Claudon}
\address{Fr\' ederic \textsc{Campana}, Beno\^it \textsc{Claudon}, Universit\' e de Lorraine, Institut \' elie Cartan Nancy, UMR 7502, B.P. 70239, 54506 Vand\oe uvre-l\`es-Nancy Cedex, France}
\email{Frederic.Campana@univ-lorraine.fr}
\email{Benoit.Claudon@univ-lorraine.fr}
\thanks{BC bénéfécie du soutien du CNRS et de l'IMPA \emph{via} un séjour de recherche effectué à l'IMPA. BC est également membre de l'ANR-10-JCJC-0111. Le présent travail a été initié alors que les auteurs se trouvaient à l'IMPA : ils souhaitent remercier l'IMPA pour son hospitalité.}
\begin{document}

\begin{abstract}
Nous montrons que les fibres du morphisme d'Albanese d'une vari\'et\'e complexe projective sp\'eciale sont sp\'eciales, r\'epondant ainsi positivement \`a une question pos\'ee dans \cite{Ca04}. L'ingr\'edient principal de la preuve est un cas particulier de la version orbifolde, conjectur\'ee dans loc.cit., de la conjecture $C_{n,m}$ d'Iitaka, \'etabli par Birkar-Chen \cite{BC13}. Nous donnons également quelques applications (groupe fondamental et revêtement universel des variétés spéciales).
\end{abstract}

\maketitle

\tableofcontents

\section{Introduction}

Nous consid\'erons ici certaines propriétés de stabilité des variétés dites `spéciales'. Rappelons que les \'el\'ements de cette classe, introduite dans \cite{Ca04}, sont les variétés projectives lisses (ou plus généralement Kählériennes compactes) qui n'admettent aucune fibration `de type général' non triviale. Une fibration $f:X\To Y$ `nette' \'etant dite de type général si la paire orbifolde $(Y,\Delta(f))$ est de type général, c'est-\`a-dire si le $\QQ$-diviseur $K_Y+\Delta(f)$ est \emph{big}. Le diviseur $\Delta(f)$ est un $\QQ$-diviseur sur la base $Y$ de la fibration $f$ encodant les fibres multiples de $f$. Pour plus de détails, nous renvoyons à l'appendice à la fin de cet article ainsi qu'aux articles \cite{Ca04} et \cite{Ca11j}.\\

Pour une classe de variétés (projectives lisses, ici) $\mathcal{S}$, nous consid\'ererons les propriétés de stabilité suivantes:
\begin{enumerate}[(i)]
\item stabilité par passage au quotient : si $f:X\To Y$ est surjective et si $X\in \mathcal{S}$, $Y$ est-elle aussi dans $\mathcal{S}$ ?
\item stabilité par passage aux fibres d'une fibration : si $f:X\To Y$ est une fibration (morphisme surjectif à fibres connexes) et si $X\in \mathcal{S}$, les fibres générales de $f$ sont-elles dans $\mathcal{S}$ ?
\item stabilité par revêtement étale : si $\pi:Y\To X$ est un revêtement étale fini et si $X\in\mathcal{S}$, $Y$ est-elle aussi dans $\mathcal{S}$ ?
\item stabilité par revêtement universel (généralisation de la question précédente) : si $X$ et $Y$ ont des revêtements universels (analytiquement) isomorphes et si $X\in\mathcal{S}$, $Y$ est-elle aussi dans $\mathcal{S}$ ?
\item stabilité par déformation : si $X\To B$ est une famille lisse et si $X_b\in\mathcal{S}$ pour un point $b\in B$, cela implique-t-il que $X_b\in\mathcal{S}$ pour tout point $b\in B$ ? Cette dernière question se scinde en deux : montrer que propriété est ouverte et fermée.
\end{enumerate}

Si $\mathcal{S}$ est la classe des variétés de type général la propriété $(i)$ n'est pas vérifiée mais les autres propriétés de stabilité $(ii)-(v)$ sont vérifiées. Le point $(iv)$ est un résultat de H. Tsuji \cite{T96}, $(v)$ \'etant l'invariance des plurigenres de Y.-T. Siu \cite{S98,S02}.

Si $\mathcal{S}$ est la classe des variétés spéciales, la propriété $(i)$ est trivialement vérifiée par définition, mais bien s\^ur pas la propriété $(ii)$ (pinceau de courbes de genre au moins 2 dans $\PP^2$). Cependant, $(ii)$  est satisfait si $f$ est l'application d'Albanese de $X\in \mathcal{S}$, par le th\'eor\`eme \ref{thm intro} ci-dessous. 

La stabilité par revêtement étale est établie dans \cite{Ca04}: la démonstration repose sur la généralisation orbifolde de la profonde faible-positivit\'e des images directes des faisceaux pluricanoniques, due \`a E.Viehweg. Le théorème \ref{thm intro} permet d'en déduire (moyennant la conjecture d'Ab\'elianit\'e sur le groupe fondamental des variétés spéciales, voir la section \ref{section abel}) une réponse affirmative à la question $(iv)$.

La stabilité par déformation, si vraie, semble pour le moment hors d'atteinte avec les techniques actuelles. Dès la dimension 3, cette question est ouverte (pour les surfaces, la stabilité par déformation du caractère spécial résulte de sa caractérisation par la dimension de Kodaira et le groupe fondamental, qui sont invariants par d\'eformation, voir \cite[Cor.3.32, p.552]{Ca04}).

\begin{thm*}\label{thm intro}
Si $X$ est une variété projective lisse et spéciale, les fibres générales de l'application d'Albanese de $X$ sont encore des variétés spéciales.
\end{thm*}

\section{Fibration d'Albanese}

\begin{thm}\label{alb surjective} 
Si $X$ est une vari\'et\'e k\"ahl\'erienne compacte sp\'eciale\footnote{L'\'enonc\'e et d\'emonstration du th\'eor\`eme \ref{tpf} sont valables pour les deux variantes, \emph{classique} et \emph{non-classique} de cette notion. Nous donnerons la d\'emonstration dans le cas non-classique. L'adaptation au cas classique se fait en rempla\c cant les termes \emph{inf} par \emph{pgcd} aux endroits ad\'equats. En revanche, notre d\'emonstration du th\'eor\`eme \ref{tfas} ne s'applique qu'aux vari\'et\'es non-classiquement sp\'eciales.}, son application d'Albanese
$$\alpha_X:X\To\alb(X)$$
est une fibration et est sans fibre multiple en codimension un.
\end{thm}

Ce r\'esultat (\'etabli dans \cite[Prop.5.3, p.576]{Ca04}) g\'en\'eralise et renforce l'\'enonc\'e similaire de Y. Kawamata \cite{K81} concernant les vari\'et\'es projectives de dimension de Kodaira nulle (qui sont sp\'eciales, par \cite[Th.5.1, p.575]{Ca04}), ainsi que celui de K. Ueno \cite[Lem.13.1, Lem.13.6]{U75} concernant les vari\'et\'es complexes compactes de dimension alg\'ebrique nulle (\'egalement sp\'eciales par \cite{Ca04}).

Dans \cite[Question 5.4, p.577]{Ca04}, il est demandé si les fibres g\'en\'eriques de l'application d'Albanese de $X$ sont sp\'eciales (pour $X$ kählérienne compacte sp\'eciale) et il est remarqu\'e que l'on peut le d\'emontrer si la conjecture $C_{m,n}^{orb}$ est vraie (cette conjecture est \'enonc\'ee en \cite[Conj.4.1, p. 564]{Ca04}).

\begin{conj}[Conjecture $C^{orb}_{n,m}$]\label{c_n,m orbifolde}
Soit $(X,D)$ une paire orbifolde lisse (compacte K\"ahler) et $f:X\to Y$ une fibration. Soit $D_y$ la restriction de $D$ \`a la fibre g\'en\'erale $X_y$ de $f$, et $(Y,D_{f,D})$ la base orbifolde de $(f,D)$ (voir annexe en fin de texte). Alors $\kappa(X,D)\geq \kappa(X_y,D_y)+\kappa(Y,D_{f,D})$ si $(Y,D_{f,D})$ est lisse, si la r\'eunion du support de $D$ et du discriminant de $f$ est \`a croisements normaux, et s'il existe une modification $u:X\to X_0$ avec $X_0$ lisse telle que tous les diviseurs $f$-exceptionnels de $X$ soient $u$-exceptionnels.
\end{conj} 

Un cas particulier de cette conjecture $C_{m,n}^{orb}$ a \'et\'e \'etablie sous une forme un peu plus faible dans \cite{BC13}:

\begin{thm}[\cite{BC13}]\label{additivit\' e orbifolde BC}
Soit $f:X\To A$ une fibration entre vari\' et\' es projectives et supposons $X$ munie d'une structure orbifolde $(X,D)$ lisse\footnote{On se ram\`ene à ce cas \`a l'aide d'une log-r\'esolution lorsque $(X,D)$ est klt, ce qui est l'hypoth\`ese de l'\'enonc\'e original de \cite{BC13}.} v\'erifiant la propri\' et\' e suivante: les fibres orbifoldes g\' en\' erales $(X_a,D_{X_a})$ sont de type g\' en\' eral. Si $A$ est de dimension d'Albanese maximale\footnote{Cela signifie que le morphisme d'Albanese de $A$ est génériquement fini sur son image.}(en particulier si $A$ est une vari\'et\'e ab\'elienne), alors:
$$\kod{X,\Delta}\ge \dimm{X}-\dimm{A}+\kod{A}.$$
\end{thm}

Nous en d\' eduisons :

\begin{thm}\label{tfas}
Si $X$ est une vari\' et\' e projective lisse et sp\' eciale, les fibres g\' en\' erales de l'application d'Albanese sont \'egalement sp\'eciales.
\end{thm}

\begin{rem}\label{cas kappa=0}
Si $X$ est une vari\' et\' e projective avec $\kappa(X)=0$, les fibres de l'application d'Albanese ont aussi une dimension de Kodaira nulle, par \cite{CH11}.
\end{rem}

Nous d\'eduirons en fait le th\'eor\`eme \ref{tfas} d'un renforcement de l'\' enonc\'e de \cite{BC13}:

\begin{thm}\label{tpf}
Soit $f:X\To A$ une fibration entre vari\' et\' es projectives et supposons $X$ munie d'une structure orbifolde lisse $(X,D)$ dont les fibres orbifoldes g\' en\' erales $(X_a,D_a)$ sont de type g\' en\' eral. Si $A$ est une vari\'et\'e de dimension d'Albanese maximale, alors : \emph{apr\`es modification birationnelle} de $(X,D)$, la base\footnote{Voir l'appendice pour les notions utilisées.} orbifolde $(Z,D_Z)$ de la fibration d'Iitaka-Moishezon $g:(X,D)\to Z$ de $(X,D)$est lisse, de type g\'en\'eral, et:
$$\dimm{Z}=\kod{X,D}\geq \dimm{X}-\dimm{A}.$$

De plus:  apr\`es un changement de base \'etale fini $A'\to A$ ad\'equat, $(X,D)$ est birationnelle au produit $(Z,D_Z)\times B'$, o\`u $B'$ est une sous-vari\'et\'e ab\'elienne de $A'$.
\end{thm}

\noindent Ce résultat est en fait une version orbifolde (plus précise) de \cite[Th. 3.1]{CCE}.

Un examen plus précis de la situation fournit aussi le renforcement suivant de \cite{BC13}, qui n'est autre que la conjecture $C^{orb}_{n,m}$ de \cite{Ca04} dans le cas particulier consid\'er\'e.
\begin{cor}\label{Cnm orbifolde}
Soit $f:(X,D)\To Y$ une fibration dont les fibres sont de type général et notons $\Delta$ la base orbifolde de $f$ induite par $D$. Si $Y$ est de dimension d'Albanese maximale, nous avons alors l'inégalité suivante entre les dimensions de Kodaira :
$$\kod{X,D}\geq \kod{X_y,D_y}+\kod{Y,\Delta}=\dim(f)+\kod{Y,\Delta}.$$
\end{cor}

\begin{ex}\label{ex} Dans la situation du th\'eor\`eme \ref{tpf}, il n'est en g\'en\'eral pas vrai que la base orbifolde $(Z,D_Z)$ est de type g\'en\'eral sans modifier $(X,D)$, même si la fibration $g$ est \emph{nette} (voir la définition \ref{defi fibration nette}). Consid\'erons $X_0=\PP^2\times E$, o\`u $E$ est une courbe elliptique, et soit $G\subset \PP^2$ une courbe lisse de degr\'e $d>3$ munie d'un morphisme surjectif $\psi:G\to E$. On munit $X_0$ du diviseur orbifolde $c.(G\times E)$, o\`u $1>c>\frac{3}{d}$ est rationnel. Donc $(Z,D_Z):=(\PP^2,c.G)$ est de type g\'en\'eral. On consid\`ere $(X_0,D_0):=(X_0, g^*(cG))$, avec les projections naturelles $g_0:X_0\to Z:=\PP^2$ et $f_0:X_0\to E$. Alors $(Z,D_Z)$ est la base orbifolde de $g_0:(X_0,D_0)\to Z$, et les fibres orbifoldes de $F:(X_0,D_0)\to E$ sont \'evidemment de type g\'en\'eral.

Soit $\eta:X\to X_0$ l'\'eclatement de $X_0$ le long de la courbe $G'\subset X_0$ graphe du morphisme $\psi: G\to E$, soit $D$ le transform\'e strict de $D_0$ dans $X$ et $f:=f_0\circ\eta:X\to E$ et $g:=g_0\circ\eta: X\to Z$ les fibrations induites sur $X$. On v\'erifie ais\'ement que les fibres orbifoldes de $f$ sont encore de type g\'en\'eral. Toutefois, dans cette situation, la base orbifolde de $g:X\to Z$ est maintenant $(Z,0)=(\PP^2,0)$ et n'est donc pas de type g\'en\'eral. Si on munit cependant le diviseur exceptionnel $E$ de $\eta$ d'une multiplicit\'e $c'$ plus grande que $\frac{3}{d}$, alors la base orbifolde de $g':(X, D+c'E)\to Z$ est $(Z, c'G)$ est de type g\'en\'eral.
\end{ex}

\begin{ex}\label{exemple bete}
L'exemple suivant montre la néc\'essité d'introduire la base orbifolde de $g:(X,D)\to Z$, m\^eme lorsque $D=0$. Considérons en effet
$$X:=\left(C\times E\right)/\langle s\times t\rangle$$
où $C$ est une courbe hyperelliptique (de genre au moins 2) d'involution $s$ et $E$ une courbe elliptique munie de $t$ une translation par un point d'ordre 2 (dans cet exemple $D=0$). La deuxième projection
$$p_2:X\To E/\langle t\rangle$$
est une fibration sur une courbe elliptique dont les fibres sont de type général et la fibration d'Iitaka est donnée par la première projection
$$p_1:X\To C/\langle s\rangle\simeq \PP^1.$$
La base de cette fibration ne devient de type général qu'une fois ajouté le diviseur orbifolde (ici le diviseur de branchement de $C\To\PP^1$).
\end{ex}

Nous d\'emontrerons les th\'eor\`emes \ref{tpf} et \ref{tfas} et le corollaire \ref{Cnm orbifolde} dans la section \ref{demo}. Avant cela, donnons quelques applications du théorème \ref{tfas}, notamment au problème envisagé dans l'introduction.

\section{Conjecture d'Ab\' elianit\' e}\label{section abel}

Une variété spéciale étant conjecturalement (par abondance) une extension successive d'orbifoldes soit Fano, soit \`a à fibré canonique numériquement trivial (en un sens birationnel ad\'equat), son groupe fondamental devrait être presque ab\'elien. La conjecture suivante, extraite de \cite{Ca04}, est bas\'ee sur cette observation.

\begin{conj}\label{conj abel} Soit $X$ une vari\' et\' e sp\' eciale. Alors:
\begin{enumerate}
\item $\pi_1(X)$ est presque ab\' elien.
\item Si $\tilde{q}(X)=0$, $\pi_1(X)$ est fini\footnote{Rappelons que l'irrégularité maximale est définie par $$\tilde{q}(X)=\mathrm{sup}\left(q(\tilde{X})\vert \tilde{X}\To X\,\textrm{étale, fini}\right).$$}.
\end{enumerate}
\end{conj}

Il est clair que la conjecture (\ref{conj abel}.1) implique la conjecture (\ref{conj abel}.2). Le théorème \ref{tfas} permet cependant de montrer que ces deux énoncés sont en fait équivalents et donc de réduire la conjecture (\ref{conj abel}.1) à l'énoncé (\ref{conj abel}.2).

\begin{prop}\label{reduction abelianite}
La validit\' e de la conjecture (\ref{conj abel}.2) implique celle de la conjecture (\ref{conj abel}.1).
\end{prop}
\begin{demo}
On proc\`ede par r\' ecurrence sur la dimension de $X$. Si $X$ v\' erifie $\tilde{q}(X)>0$, on choisit un rev\^etement \' etale (toujours not\' e $X$) pour lequel $q(X)=\tilde{q}(X)$. C'est possible, puisque les rev\^etements \'etales de $X$ restent sp\'eciaux et ont une irr\'egularit\'e major\'ee par $\dim(X)$. On consid\`ere l'application d'Albanese de $X$. Comme $\alpha_X$ est sans fibre multiple en codimension 1, la suite de groupes:
$$\pi_1(F)\To \pi_1(X)\To \pi_1(\mathrm{Alb}(X))\To 1$$
est exacte, $F$ \'etant une fibre lisse de $\alpha_X$. D'après le théorème \ref{tfas}, $F$ est une variété spéciale. Par r\' ecurrence sur la dimension, son groupe fondamental $\pi_1(F)$ est presque ab\' elien et $\pi_1(F)_X$ (son image dans celui de $X$) l'est \' egalement. Comme il en est de m\^eme pour $\pi_1(\mathrm{Alb}(X))$, les r\' esultats de \cite{ca98} montrent que $\pi_1(X)$ est bien presque ab\' elien.
\end{demo}

Notons ici que la conjecture (\ref{conj abel}.1) est connue en dimension 3 : une variété Kählérienne compacte \emph{spéciale} de dimension 3 a un groupe fondamental presque abélien. Il s'agit en effet du résultat principal de \cite{CC14}.

\section{Rev\^etement universel}

Nous consid\'erons ici la stabilité par revêtement universel (point $(iv)$ de l'introduction): \^etre une vari\' et\' e sp\' eciale ne devrait d\' ependre que du rev\^etement universel. En effet, il est conjectur\'e \cite{Ca04} qu'\^etre spéciale \'equivaut à l'annulation de la pseudo-métrique de Kobayashi (pour une discussion de ces aspects, voir \cite{Ca04} et la question ). Cette pseudo-distance étant invariante par revêtement \cite{livreKob}, son annulation \'equivaut \`a l'annulation de celle de son revêtement universel. Ainsi, si $X$ et $Y$ sont deux variétés K\"ahl\'eriennes compactes  ayant le même revêtement universel, elles sont, selon cette conjecture, simultan\'ement sp\'eciales ou non-sp\'eciales.  

\begin{qt}\label{question rev univ}
Soient $X$ et $Y$ sont deux variétés k\"ahl\'eriennes compactes ayant le même revêtement universel. Si $X$ est spéciale, $Y$ l'est-elle également ?
\end{qt}
Cette question est similaire \`a la conjecture d'Iitaka: si $X$ est une variété abélienne, $Y$ est-elle une variété abélienne (à revêtement étale fini près)?

H. Tsuji \cite{T96} r\'esoud cette question dans le cas \og oppos\'e\fg des vari\'et\'es de type g\'eneral: il y montre en effet que, dans la situation de la question \ref{question rev univ}, si $X$ est de type général, alors $Y$ l'est aussi.\\

Le résultat suivant montre que la conjecture d'abélianité \ref{conj abel} implique une réponse positive (dans le cas projectif) à la question \ref{question rev univ}. Notons qu'il est \'etabli dans \cite{Ca04} que la conjecture d'Ab\'elianit\'e implique une r\'eponse positive \`a la conjecture d'Iitaka pr\'ec\'edente (sous une hypoth\`ese plus faible: la lin\'earit\'e du groupe fondamental de $X$, au lieu de son ab\'elianit\'e).

\begin{thm}\label{stab rev universel}
Soit $X$ et $Y$ deux vari\' et\' es projectives lisses ayant le m\^eme rev\^etement universel. Si $X$ est une vari\' et\' e sp\' eciale et si son groupe fondamental est presque ab\' elien, alors $Y$ est sp\' eciale \' egalement.
\end{thm}
En combinant le résultat ci-dessus avec le résultat principal de \cite{CC14} (valable dans le cas K\"ahl\'erien compact), nous obtenons :
\begin{cor}\label{Corollaire rev univ}
Si $X$ et $Y$ sont deux variétés Kählériennes compactes ayant le même revêtement universel et de dimension au plus 3, alors $X$ est spéciale si et seulement si $Y$ l'est.
\end{cor}

\begin{demo}
Notons $\mathcal{D}$ le rev\^etement universel de $X$ et $Y$. Comme $X$ est sp\' eciale et de groupe fondamental presque ab\' elien, on peut supposer que l'application d'Albanese de $X$ induit un isomorphisme sur les groupes fondamentaux. En particulier, si $A$ d\' esigne la vari\' et\' e d'Albanese de $X$, l'application induite entre les rev\^etement universels
$$\xymatrix{\mathcal{D}\ar[d]_{\pi_X}\ar[r]^{\tilde{\alpha}} & \tilde{A}\simeq \CC^{q(X)}\ar[d]\\
X\ar[r]^{\alpha} & A
}$$
est propre. Nous savons en outre que les fibres g\' en\' erales de $\tilde{\alpha}_X$ sont sp\' eciales (puisque sont celles de $\alpha$). Comme $\tilde{\alpha}$ est propre (et $\tilde{A}$ \' etant de Stein), les automorphismes de $\mathcal{D}$ induisent\footnote{Le translaté d'une fibre par un automorphisme de $\mathcal{D}$ est encore une sous-vari\' et\' e compacte connexe et son image par $\tilde{\alpha}_X$ est donc un point.} des automorphismes de $\tilde{A}$ ; en particulier, $\tilde{\alpha}$ est donc équivariant pour une action (proprement discontinue) de $\pi_1(Y)$ sur $\tilde{A}$. En passant au quotient, nous obtenons ainsi un morphisme :
$$f:\mathcal{D}/\pi_1(Y)=Y\To Z:=\tilde{A}/\pi_1(Y).$$

Raisonnons par l'absurde et supposons que $Y$ n'est pas sp\'eciale. Soit alors
$$\cc:Y\To (C,\Delta_C)$$
un repr\'esentant net de son c{\oe}ur, l'orbifolde $(C,\Delta_C)$ \'etant lisse (c'est-\`a-dire que $C$ est lisse et le support de $\Delta_C$ \`a croisements normaux) et de type g\'en\'eral avec $p:=\dim(C)>0$. En particulier, les formes m\'eromorphes ayant des pôles le long de $\Delta_C$ deviennent donc holomorphes sur $Y$ (proposition \ref{relevement base orbi nette}) :
$$\cc^*\left(H^0(C,m(K_C+\Delta_C))\right)\subset H^0(Y,(\Omega^p_Y)^{\otimes m}).$$
En prenant l'image réciproque par $\pi_Y$, il en résulte que les pluriformes méromorphes sur $C$ se relèvent en des tenseurs holomorphes sur $\mathcal{D}$ le revêtement universel de $Y$ :
$$(\cc\circ \pi_Y)^*\left(H^0(C,m(K_C+\Delta_C))\right)\subset H^0\left(\mathcal{D},(\Omega^p_{\mathcal{D}})^{\otimes m}\right).$$

Les fibres de $f$ correspondant par projection \`a celles de $\tilde{\alpha}$,  elles sont sp\' eciales. Par propriété universelle du c{\oe}ur, il existe une factorisation $g:Z\dasharrow C$ de $c$ telle que $\cc=g\circ f$. Résumons la situations dans le diagramme suivant :
\begin{equation}\label{diagramme rev universel}
\begin{gathered}
\xymatrix{\mathcal{D}\ar[rr]^{\tilde{\alpha}}\ar[d]_{\pi_Y} && \tilde{A}\ar[d]^{\rr}\\
Y\ar[rr]^{f}\ar[rd]_{\cc} & & Z\ar[ld]^{g}\\
&C&
}
\end{gathered}
\end{equation}
o\`u $\rr:\tilde{A}\to Z$ et $\pi_Y:\mathcal{D}\to Y$ sont les quotients par $\pi_1(Y)$. Nous avons donc : $\cc\circ \pi_Y=g\circ \rr\circ\tilde{\alpha}$ et il s'ensuit que :
$$(\cc\circ \pi_Y)^*(H^0(C,m(K_C+\Delta_C))=\tilde{\alpha}^*((g\circ \rr)^*\left(H^0(C,m(K_C+\Delta_C))\right)\subset H^0(\mathcal{D},(\Omega^p_{\mathcal{D}})^{\otimes m}).$$
Dans notre situation, les pluriformes méromorphes se relèvent en fait déjà comme des tenseurs holomorphes sur $\tilde{A}$ :
$$(g\circ \rr)^*\left(H^0(C,m(K_C+\Delta_C))\right)\subset H^0(\tilde{A},(\Omega^p_{\tilde{A}})^{\otimes m}).$$
En effet, le morphisme $\tilde{\alpha}:\mathcal{D}\to \tilde{A}$ n'ayant pas de fibres multiples en codimension $1$ (par le th\'eor\`eme \ref{alb surjective}), les formes différentielles m\'eromorphes $(g\circ \rr)^*(m(K_C+\Delta_C))$ se prolongent en codimension $1$ (et donc partout, par le th\'eor\`eme d'Hartogs) en des sections holomorphes de $(\Omega^p_{\tilde{A}})^{\otimes m}$. Pour se convaincre de cela, considérons $a\in \tilde{A}$ un point où $(g\circ \rr)$ est définie et qui a pour image un point du support de $\Delta_C\subset C$, et $s$ une forme méromorphe section de $(g\circ \rr)^*(m(K_C+\Delta_C))$. Comme nous l'avons constaté ci-dessus, le tenseur $\tilde{\alpha}^*(s)$ devient holomorphe sur $\mathcal{D}$. On peut cependant, par absence de fibres multiples en codimension $1$ de $\tilde{\alpha}$, supposer l'existence d'une section locale\footnote{La consid\'eration des multiplicit\'es non classiques est ici essentielle.} $\sigma: \tilde{A}\to \mathcal{D}$ de $\tilde{\alpha}$ qui permet de redescendre $\tilde{\alpha}^*(s)$ sur $\tilde{A}$ en appliquant $\sigma^*$.\\
La propri\'et\'e
$$(g\circ \rr)^*\left(H^0(C,m(K_C+\Delta_C))\right)\subset H^0(\tilde{A},(\Omega^p_{\tilde{A}})^{\otimes m})$$
\'etant \'etablie, nous pouvons appliquer la proposition \ref{KO orbifolde} ci-dessous (qui n'est autre que \cite[Th.8.2, p.599]{Ca04}), en y choisissant $U=\tilde{A}=V-D$, $V=\PP^q, D=\PP^{q-1}$ : les pluriformes $(g\circ \rr)^*\left(H^0(C,m(K_C+\Delta_C))\right)$ se prolongent \`a $V$ en des \'el\'ements de
$$H^0(V,(\Omega^p_{V})^{\otimes m}\otimes \mathcal O_V((m-1)D)).$$
Cependant, l'espace en question est nul si $V=\PP^q$ et $D=\PP^{q-1}$ \cite[Ex.8.8, 8.13 et Prop.8.14]{Ca04}. Ceci contredit le fait que $(C,\Delta_C)$ est de type g\'en\'eral avec $p>0$.
\end{demo}

Nous avons utilisé ci-dessus le fait que l'on peut étendre les sections et contrôler l'ordre des pôles pour une application non dégénérée vers une orbifolde de type général, énoncé qui est une version orbifolde du théorème de Kobayashi-Ochiai \cite{KO75}. Il est à noter qu'une version orbifolde (plus faible, et insuffisante pour conclure ici) avait déjà été établie dans \cite{Sak}.

\begin{prop}\label{KO orbifolde}
Soit $f:U=V\backslash D\dashrightarrow (Z,\Delta)$ une application méromorphe non-dégénérée vers une orbifolde lisse de type général avec $D$ un diviseur de la variété lisse (non n\'ecessairement compacte) $V$. Si les formes méromorphes (pour $m>0$ suffisamment divisible)
$$H^0(Z,m(K_Z+\Delta))$$
se relèvent par $g^*$ en des tenseurs \textrm{holomorphes} sur $U$, alors on a en fait :
$$g^*H^0(Z,m(K_Z+\Delta))\subset H^0\left(V,(\Omega^p_{V})^{\otimes m}\otimes \mathcal O_V((m-1)D)\right).$$
\end{prop}
Les arguments donnés dans la section 8 de \cite{Ca04} fournissent en effet cet énoncé. Rappelons en les grandes lignes. Le fait de supposer l'orbifolde $(Z,\Delta)$ de type général permet de construire une pseudo-forme volume dont la courbure de Ricci est négative. Cette pseudo-forme volume étant construite à partir d'une base de sections de $m(K_Z+\Delta)$ (pour $m$ assez grand et suffisamment divisible), l'hypothèse faite dans l'énoncé ci-dessus montre que l'image réciproque par $f$ de cette pseudo-forme volume a un sens et est encore une pseudo-forme volume à courbure de Ricci négative (au moins dans le cas équidimensionnel, cas qui nous occupe ici). On applique alors le lemme d'Ahlfors-Schwarz dans un polydisque contenu dans $V$ et centré sur une composante de $D$ : la pseudo-forme volume est dominée par la forme volume standard du polydisque (épointé) et en particulier intégrable au voisinage de $D$. Ceci montre immédiatement que les sections s'étendent et un calcul local \cite[Prop.8.28]{Ca04} permet de contrôler l'ordre des pôles.

\section{D\'emonstrations} \label{demo}

\noindent Nous commençons par donner la\\

\begin{demo}[du théorème \ref{tpf}]
Tout d'abord remarquons que l'on peut remplacer $A$ par sa variété d'Albanese et donc supposer que $A$ est une variété abélienne. En effet, cela altère seulement la connexité des fibres générales de $f$ mais n'est pas utilis\'e dans la démonstration.\\

\noindent\textbf{Première étape : étude de la fibration d'Iitaka de \mathversion{bold}$(X,D)$\mathversion{normal}.} Comme la paire orbifolde initiale est lisse, nous pouvons la remplacer par un modèle lisse birationnel au sens de la définition \ref{defi morphisme orbi} et
qui sera encore noté $(X,D)$ : cela ne change pas les hypothèses de départ et la dimension de Kodaira $\kappa(X,D)$ est également inchangée car la paire initiale est lisse. Nous supposerons donc dorénavant que le morphisme d'Iitaka $g:(X,D)\to Z$ de $K_X+D$ est holomorphe et \emph{net} au sens de \cite{Ca11j}. Ce morphisme est non trivial puisque le théorème \ref{additivit\' e orbifolde BC} appliqué à $f$ fournit la minoration
$$\kappa(X,D)\geq \dim(f)+\kappa(A)>0.$$
Nous munissons alors la base $Z$ de $g$ de la structure orbifolde $(Z,D_Z)$ induite par $g$ et $D$ sur $X$ (comme dans la définition \ref{base orbi}). 

Soit $(X_z,D_z:=D\cap X_z)$ une fibre orbifolde g\'en\'erique de $g$ (lisse par le th\'eor\`eme de Sard). Par définition,  nous avons $\kappa(X_z,D_z)=0$ et $(X_z,D_z)$ est donc sp\'eciale \cite{Ca11j}, ainsi {\it a fortiori} que $X_z$. Il s'ensuit que $f(X_z)=C_z\subset A$ est une sous-vari\'et\'e sp\'eciale de $A$. C'est donc une sous-vari\'et\'e ab\'elienne de $A$, constante \`a translation pr\`es. Soit $q:A\to B:=A/C$ la vari\'et\'e ab\'elienne quotient. Il existe donc un morphisme $h:Z\to B$ tel que $q\circ f=h\circ g:X\to B$.

Nous allons montrer, c'est l'\'etape essentielle, que l'application
$$f_z:=f_{\vert X_z}:X_z\to C_z$$
est g\'en\'eriquement finie (pour $z\in Z$ g\'en\'eral). Pour $a\in C_z$ g\'en\'eral, notons $(X_{z,a},D_{z,a})$ la fibre orbifolde (lisse) de $f_z:(X_z,D_z)\to C_z$. La fibre orbifolde g\'en\'erique $(X_a,D_a)$ de $f$ est par hypothèse de type général et est recouverte par les $(X_{z,a},D_{z,a})$ pour les $z$ tels que $a\in C_z$, c'est-à-dire tels que $h(z)=q(a)$, ou encore: $z\in h^{-1}(q(a))$. Nous en déduisons que $(X_{z,a},D_{z,a})$ est encore de type général \cite[Th.9.12, p.893]{Ca11j} et nous pouvons donc appliquer de nouveau le théorème \ref{additivit\' e orbifolde BC} à $f_z$. Ceci donne alors $$0=\kappa(X_z,D_z)\geq \kappa(X_{z,a},D_{z,a})+\kappa(C_z)=\dimm{X_{z,a}}\geq 0$$
et implique que $\dimm{X_{z,a}}=0$, c'est-à-dire que $f_z:X_z\to C_z$ est g\'en\'eriquement finie.

Puisque $\kappa(X_z,D_z)=0$ et que $C_z$ est une variété abélienne, on en d\'eduit que la factorisation de Stein de $f_z$
$$f_z:X_z\stackrel{f'_z}{\To}C'_z\stackrel{p_z}{\To}C_z$$
vérifie : $p_z$ est \'etale et $D_z$ est $f'_z$-exceptionnel. L'application naturelle
$$\eta:=g\times f:X\to Z\times_B A$$ est donc surjective et g\'en\'eriquement finie. De plus, $g(D)$ est contenu dans un diviseur de $Z$.\\

\noindent\textbf{Deuxième étape : cas où $\eta$ est birationnelle.} Supposons tout d'abord que $\eta$ est birationnelle et voyons comment conclure dans ce cas.

Nous munissons alors $X_0:=Z\times_B A$ du diviseur $\eta_*(D):=D_0$ et examinons la projection naturelle $g_0:X_0\to Z$. Puisque $g_0$ est submersive et que $g_0(D_0)=g(D)$ est un diviseur de $Z$, nous avons
$$D_0=g_0^*(D^0_Z)$$
si $D^0_Z$ est la base orbifolde de $g_0:(X_0,D_0)\to Z$. De plus,
$$K_{X_0}+D_0=(g_0)^*(K_{Z_0}+D_Z^0)$$
et les dimensions de Kodaira coïncident donc $\kappa(X_0,D_0)=\kappa(Z,D^0_Z)$. Puisque $$K_{X}+D=\eta^*(K_{X_{0}}+D_{0})+E,$$
o\`u $E$ est diviseur $\eta$-exceptionnel (pas n\'ecessairement effectif) et que $X_{0}$ est lisse, le th\'eor\`eme de Hartogs montre que
$$\dimm{Z}:=\kappa(X,D)\leq \kappa(X_0,D_0)=\kappa(Z,D^0_Z)$$
et $(Z,D^0_Z)$ est donc de type g\'en\'eral.\\

\noindent\textbf{Troisième étape : changement de base.} Nous allons montrer que $\eta$ devient birationnelle apr\`es changement de base \'etale fini ad\'equat $\pi:A'\to A$ et analyserons ce changement de base.

Soit, en effet, $f_*(\pi_1(X_z))\subset \pi_1(C_z)$ : c'est un sous-groupe d'indice fini de $\pi_1(C_z)=\pi_1(C)$ ind\'ependant de $z$ g\'en\'erique dans $Z$. Il existe un rev\^etement \'etale fini $\pi: A'\to A$ tel que $\pi_1(A')\cap \pi_1(C_z)=F_*(\pi_1(X_z))=\pi_1(C')$ pour une sous-vari\'et\'e ab\'elienne $C'$ de $A'$ et tel que, de plus, $A'/C'=A/C=B$. Ceci est une conséquence du fait que $\pi_1(A)$ est (bien que non-canoniquement) isomorphe \`a $\pi_1(C)\times \pi_1(B)$.

Nous considérons alors la fibration
$$f':\left(X',D':=\pi^*(D)\right)\To A'$$
obtenue par changement de base $A'\To A$ et notons $\pi$ le revêtement étale $X'\To X$. Une fois ce changement de base effectué, l'application
$$\eta':=g'\times f' :X'\To Z'\times_{B}A'=:X'_0$$
est alors surjective et envoie birationnellement $X'_z$ sur $C'_z$ et $\eta'$ est donc birationnelle.  Remarquons toutefois que, lors de cette opération, le morphisme induit (avec les notations évidentes) $Z'\To Z$ n'est en général pas étale.

Consid\'erons enfin la factorisation de Stein
$$\eta:X\stackrel{\varphi}{\To} X_1\stackrel{p}{\To}X_0$$
de $\eta=p\circ \varphi$. L'application $\eta$ \'etant génériquement finie, $\varphi$ est birationnelle (connexe) et $p$ finie. Nous allons voir plus loin que $p$ est galoisienne au-dessus de $X_0$ pour une action ad\'equate naturelle de $G:=\mathrm{Gal}(C'/C)=\mathrm{Gal}(A'/A)$ sur $X_1$.
Nous pouvons résumer la situation dans le digramme suivant.
\begin{equation}\label{diagramme chgt base}
\begin{gathered}
\xymatrix{X'\ar[rrr]^{\pi}\ar[ddr]^{\eta'}\ar[ddd]_{g'} & & & X\ar[ld]_{\varphi}\ar[ldd]^{\eta}\ar[ddd]^{g}\\
&&X_1\ar[d]_{p}&\\
&X'_0\ar[ru]^q\ar[r]^{\pi_0}\ar[ld]^{g'_0}&X_0\ar[rd]_{g_0}&\\
Z'\ar[rrr]^{\rho} & & & Z.
}
\end{gathered}
\end{equation}
Les applications qui apparaissent ci-dessus poss\'edent les propriétés suivantes:
\begin{itemize}\label{liste}
\item $g$, $g'$, $g_0$ et $g'_0$ sont des fibrations, les deux dernières sont de plus lisses.
\item $\pi$, $\pi_0$ et $\rho$ sont finies (galoisiennes), $\pi$ étant étale.
\item $\eta'$ et $\varphi$ sont birationnelles.
\item $p$ est finie.
\end{itemize}

Le groupe $G$ agit de deux façons.  Tout d'abord, l'action de $G$ sur $X'$ (déduite par changement de base de celle de $G$ par translations sur $A'$) induit une action de $G$ sur $Z'$ qui est la base de la fibration d'Iitaka du $\QQ$-diviseur $G$-invariant $K_{X'}+D'$. Il est à noter que cette action sur $Z'$ possède en général des points fixes. D'autre part, il est aisé de constater que $G\times G$ agit sur $X'_0$. En effet, la variété en question est obtenue comme produit fibré :
$$X'_0=Z'\times_B A'$$
et le groupe $G$ agit à la fois sur $Z'$ et $A'$. Comme les projections $Z'\To B$ et $A'\To B$ sont invariantes sous ces actions, on en déduit une action compatible de $G\times G$ sur $X'_0$. Notons au passage que le quotient de $X'_0$ par l'action de $G\times G$ n'est autre que $X_0$. Nous restreignons désormais cette action à la diagonale de $G\times G$. C'est en effet cette action qui est naturelle dans le cadre qui est le nôtre, comme le montre l'affirmation ci-dessous.
\begin{fact*}
Les applications $\pi$ et $g'$ sont $G$-équivariantes (pour les actions pr\'ec\'edentes de $G$ sur $X'$, $A'$ et $Z'$). En particulier, l'application $\eta'$ est $G$-équivariante pour l'action naturelle de $G$ sur $X'$ et pour l'action diagonale de $G$ sur $X'_0$.
\end{fact*}

Remarquons enfin qu'il existe un morphisme naturel $q:X'_0\to X_1$. Ceci r\'esulte du caractère fonctoriel de la factorisation de Stein. En effet, l'application $\eta'$ est $G$-équivariante d'après l'affirmation ci-dessus et nous pouvons passer au quotient. Nous obtenons ainsi une application birationnelle
$$\eta'/G:X'/G=X\To X'_0/G$$
qui n'est autre que l'application $\varphi$ par unicité de la factorisation de Stein. Nous en déduisons en particulier que $X_1$ n'est autre que le quotient\footnote{Remarquons que, le groupe $G$ étant abélien, $G$ est distingué dans $G\times G$ et on en déduit que $X_1$ hérite \emph{via} $X'_0$ d'une action de $G\times G/G\simeq G$. Cette action peut également être décrite de la manière suivante. Soit $Z^*\subset Z$ l'ouvert de Zariski dense (puisque $X_z\to A_b,b:=h(z)$ est \'etale pour $z\in Z$ g\'en\'erique) au-dessus duquel $g:X\to Z$ est lisse et soit $X_1^*:=g^{-1}(Z^*)$. L'application $g_1^*:X_1^*\to Z^*$ est donc lisse : c'est un fibr\'e principal de fibre $C'$ au-dessus de $Z^*$. Le groupe $G$ agit sur $X_1^*$ par translations dans les fibres de $g_1^*$. Cette action se prolonge analytiquement \`a $X_1$. En effet, si nous notons $\Gamma^*_{\gamma},\,\gamma\in G$ le graphe dans $X^*_1\times X_1^*$ de l'automorphisme de $X_1^*$ défini ci-dessus, il est contenu dans $(p\times p)^{-1}(\Delta)$ (avec $\Delta$ la diagonale de $X_0\times X_0$), l'action de $G$ laissant invariante la projection $p$. Son adh\'erence dans $X_1\times X_1$ est donc algébrique, puisque $p\times p:X_1\times X_1\to X_0\times X_0$ est finie. C'est donc le graphe d'un automorphisme de $X_1$, puisque $X_1$ est normal.} de $X'_0$ sous l'action diagonale de $G$. Comme l'action sur le deuxième facteur (celle qui vient de $A'$) est sans point fixe, on en déduit que $q:X'_0\To X_1$ est étale et que $g_1:X_1\To Z$ est obtenu en prenant le quotient de $g'_0:X'_0\To Z'$.\\

\noindent\textbf{Conclusion.} Considérons alors les orbifoldes $(X'_0,D'_0:=\eta'_*(D'))$ et $(X_1,D_1:=\varphi_*(D))$ et notons que le morphisme
$$(X'_0,D'_0)\stackrel{q}{\To}(X_1,D_1)$$
est étale au sens orbifolde, \emph{i.e.}
$$K_{X'_0}+D'_0=q^*\left(K_{X_1}+D_1\right)$$
car les composantes contractées par $\varphi$ sur $X$ et par $\eta'$ sur $X'$ se correspondent \emph{via} $\pi$.

Munissons alors $Z'$ et $Z$ des bases orbifoldes de $g'_0:(X'_0,D'_0)\To Z'$ et $g_1:(X_1,D_1)\To Z$ ; nous les noterons $(Z',D^0_{Z'})$ et $(Z,D^1_Z)$. Il nous suffit alors de vérifier que
$$\rho:(Z',D^0_{Z'})\To (Z,D^1_Z)$$
est étale (au sens orbifolde pr\'ec\'edent) en codimension un. En effet, une fois ceci établi, nous obtiendrons l'égalité souhaitée des dimensions de Kodaira :
$$\dim(Z)=\kod{X,D}\leq\kod{X_1,D_1}=\kod {X'_0,D'_0}=\kod{Z',D^0_{Z'}}\leq\dim(Z).$$

Un calcul simple de multiplicités montre l'égalité :
$$K_{Z'}+D^0_{Z'}=\rho^*\left(K_Z+D^1_Z\right).$$
En effet, nous avons vu ci-dessus que la fibration $g_1:X_1\To Z$ est obtenue à partir de la fibration lisse $g'_0:X'_0\To Z'$ en prenant le quotient par l'action de $G$. Cela montre que, si $F$ est un diviseur de $Z$, la multiplicité de $F$ dans $D^1_Z$ est donnée par :
$$m_{D^1_Z}(F)=r_F\cdot m_{D^0_{Z'}}(F')$$
où $r_F$ est l'ordre de ramification de $\rho$ le long de $F$ et $F'$ est une composante irréductible de $\rho^{-1}(F)$.

En écrivant $D^1_Z$ sous la forme
$$D^1_Z=\sum_{i\in I} (1-\frac{1}{r_im_i})F_i$$
où $(F_i)_{i\in I}$ est l'ensemble des composantes de $D^1_Z$ et $m_i$ les multiplicités correspondantes sur $Z'$, il vient alors :
\begin{align*}
\rho^*\left(K_Z+D^1_Z\right)&=\rho^*\left(K_Z+\sum_{i\in I} (1-\frac{1}{r_im_i})F_i\right)\\
&=\rho^*(K_Z)+\sum_{i\in I} r_i(1-\frac{1}{r_im_i})\rho^*(F_i) \\
&=\rho^*(K_Z)+\sum_{i\in I} (r_i-1)\rho^*(F_i)+\sum_{i\in I}(1-\frac{1}{m_i})\rho^*(F_i)\\
&=K_{Z'}+D^0_{Z'}
\end{align*}
et le morphisme $\rho$ est bien orbifolde-étale en codimension un.

La derni\`ere assertion de l'\'enonc\'e se d\'eduit du th\'eor\`eme de compl\`ete r\'eductibilit\'e de Poincar\'e en choisissant $A'=C'\times B'$. 
\end{demo}

\begin{rem}\label{rmq necessite revetement}
L'exemple \ref{exemple bete} montre qu'il est nécessaire d'effectuer un changement de base pour obtenir le caractère birationnel de l'application $\eta$.
\end{rem}

\begin{rem}
L'exemple \ref{ex} ci-dessus montre que, dans la situation pr\'esente, $(Z,D^0_Z)$ peut \^etre de type g\'en\'eral tandis que $(Z,D_Z)$ ne l'est pas, si les multiplicit\'es sur le diviseur exceptionnel de $\eta$ sont insuffisamment grandes.
\end{rem}

\noindent Nous d\'emontrons maintenant le th\'eor\`eme \ref{tfas}.\\

\begin{demo}[du théorème \ref{tfas}]
Soit $\alpha_X:X\to A:=\mathrm{Alb}(X)$ le morphisme d'Albanese de $X$ qui est supposée lisse, connexe, projective et sp\'eciale. Donc (théorème \ref{alb surjective}) $\alpha_X$ est surjective et \`a fibres connexes. Supposons que ses fibres ne soient pas sp\'eciales. Il existe alors un c{\oe}ur relatif
$$c:=c_{X/A}:X\To Y:=C(X/A)$$
et une fibration $f:Y\To A$ tels que le diagramme
$$\xymatrix{X\ar[rr]^{c}\ar[rd]_{\alpha_X} && (Y,D)\ar[ld]^f\\ &A&
}$$
commute. Le morphisme $c$ peut, quitte \`a modifier $X$, \^etre rendu holomorphe et la base orbifolde $(Y,D)$ de $c$ \^etre suppos\'ee lisse. Les fibres orbifoldes g\'en\'erales de $f:(Y,D)\to A$ sont les parties $f$-verticales de la base orbifolde de $c$, par le lemme \ref{lemme coeur relatif} ci-dessous, et sont donc de type g\'en\'eral.

Nous pouvons donc appliquer le th\'eor\`eme \ref{tpf} au morphisme $f:(Y,D)\To A$. Nous en déduisons un morphisme birationnel $\eta: Y\to Y_0$ construit comme dans la démonstration du th\'eor\`eme \ref{tpf}. En particulier, si $D_0=\eta_*(D)$, la fibration d'Iitaka-Moishezon $g_0:(Y_0,D_0)\to (Z,D^0_Z)$ est de type général et nous noterons $g:=g_0\circ\eta$ la compos\'ee.

Quitte \`a modifier $X,Y$ et $Z$, nous pouvons supposer que l'application $g\circ c:X\to Z$ est nette. Même si toute fibration admet un modèle net (voir la proposition \ref{existence modele net}), la situation étudiée présente une difficulté supplémentaire : la fibration $g$ dépend de $c$ et il n'est \emph{a priori} pas évident qu'une modification de $Z$ soit compatible avec le fait d'être la base de la fibration d'Iitaka de $(Y,D)$. Pour parvenir à nos fins (rendre $g\circ c$ net), nous pouvons par exemple choisir un mod\`ele birationnel de $c:X\to (Y,D)$ de telle sorte que $\kappa(Y,D)$ soit minimal (et rappelons que $\dim(Z)=\kappa(Y,D)>0$ d'après le théorème \ref{additivit\' e orbifolde BC}). Cette dimension de Kodaira sera pr\'eserv\'ee par toute modification nette de ce mod\`ele de $c$ (voir \cite[\S 1.3]{Ca04}). Une fois ce choix fait, nous remarquons que si $Z'\to Z$ est une modification de $Z$ et si nous considérons les changements de base qui s'en déduisent :
$$\xymatrix{X\ar[r]^{c} & Y\ar[r]^{g} & Z\\
X'\ar[r]^{c'}\ar[u] & Y'\ar[r]^{g'}\ar[u] & Z'\ar[u],
}$$
l'application $g':Y'\to Z'$ s'identifie à nouveau avec la fibration d'Iitaka de la paire $(Y',D'=\Delta(c'))$ car $\kappa(Y',D')=\kappa(Y,D)$. Fort de cette remarque, nous pouvons donc construire un modèle net de $g\circ c:X\to Z$ en commençant par aplatir cette fibration puis en désingularisant l'espace total (nous renvoyons à nouveau à \cite[Cor.1.13]{Ca04} pour les détails).
%Les morphismes $\eta: Y\to Y_0$ et $g_0:(Y_0,D_0)\to (Z,D^0_Z)$ correspondants seront donc obtenus par changements de base birationnels $(Z',D^0_{Z'})\to (Z_0,D_Z^0)$. Nous pouvons maintenant choisir ce mod\`ele $c':X'\to (Y',D')$ de telle sorte que $\eta':(Y',D')\to (Y'_0,D'_0)$ soit nette (en aplatissant le morphisme $\eta:Y\to Y_0$ initial, et en d\'esingularisant le produit fibr\'e normalis\'e). La base orbifolde $(Z',D^0_{Z'})$ de $g'_0:(Y'_0,D'_0)\to (Z'_0,D'^0_{Z'})$ correspondante sera encore de type g\'en\'eral, avec $\eta'\circ c':X'\to (Y'_0,D'_0)$ nette, ainsi donc que $g'\circ c':X'\to (Z'_0,D^0_{Z'})$.

En particulier, $Z$ est lisse et le morphisme $g_0:Y_0\to Z$ est submersif \`a fibres connexes entre variétés lisses. De plus, $(Z,D^0_Z)$ est aussi la base orbifolde de $g\circ c$, puisque, pour toute composante $F\subset Z$ du support de $D^0_Z$, et toute composante irr\'eductible $E$ de $(g\circ c)^{-1}(F)$, la multiplicit\'e de $E$ dans $(g\circ c)^*(F)$ est un multiple de celle de $E_0=g_0^{-1}(F)$ dans $g_0^*(F)$. Or cette multiplicit\'e est aussi celle de $E_0$ dans $D_0=\eta_*(D)$, puisque les fibres de $g_0$ sont toutes irréductibles.

Il en r\'esulte que $X$ n'est pas sp\'eciale car $(Z,D^0_Z)$ est de type général. Ceci constitue une contradiction et montre que les fibres g\'en\'erales de $\alpha$ sont donc bien sp\'eciales. \end{demo}

\noindent Dans la démonstration précédente, nous avons utilisé le résultat élémentaire suivant.

\begin{lem}\label{lemme coeur relatif}
Soit $f:X/A\To Y/A$ une fibration au dessus de $A$ : $X$ et $Y$ sont elles-mêmes fibrées sur $A$ et le diagramme
$$\xymatrix{X\ar[rr]^{f}\ar[rd] && Y\ar[ld]\\ &A&
}$$
commute. Considérons $D:=D_f$ le diviseur orbifolde de $f$. Si $a\in A$ est un point suffisamment général, notons $f_a:X_a\To Y_a$ la restriction de $f$ aux fibres $X_a$ et $Y_a$ et $D_a$ la restriction de $D$ à $Y_a$. Nous avons alors
$$D_a\geq D_{f_a},$$
où $D_{f_a}$ est le diviseur orbifolde de $f_a$. En particulier, si $f$ est le c{\oe}ur relatif de $X\To A$, alors les fibres de $(Y,D)\To A$ sont de type général.
\end{lem}
\begin{demo}
Soit $D_i$ une composante de $D$ et écrivons la décomposition :
$$f^*(D_i)=\sum_{j\in J} m_j^{(i)}D_j^{(i)}+R_i$$
où les diviseurs $D_j^{(i)}$ (contenus dans $X$) s'envoient surjectivement sur $D_i$ et où $R$ est $f$-exceptionnel. En restreignant cette égalité à $X_a$, on constate que l'ensemble des diviseurs sur lequel on prend l'infimum pour calculer la multiplicité de $D_{i\vert Y_a}$ dans $D_{f_a}$ est plus grand que $J$ : en effet, certaines composantes de $R_{\vert Y_a}$ peuvent ne plus être $f_a$ exceptionnelles. Cela fournit l'inégalité annoncée.
\end{demo}

Pour finir, donnons quelques mots d'explication sur la\\
\begin{demo}[du Corollaire \ref{Cnm orbifolde}]
Reprenons les notations de la démonstration du théorème \ref{tpf} : nous avons un diagramme
$$\xymatrix{(X,D)\ar[r]^{f}\ar[d]_g & (Y,\Delta)\ar[d]\\
Z\ar[r]^h & B
}$$
dans lequel $Z$ est la base de la fibration d'Iitaka associée à $K_X+D$ et $B$ un certain quotient de $Y$. Il vient donc :
$$\kod{X,D}=\dim(Z)=\dim(f)+\dim(B).$$
Pour montrer l'inégalité souhaitée, il suffit de montrer que $\dim(B)\geq\kod{Y,\Delta}$ ce qui revient à montrer que, si $b\in B$ est un point général, la dimension de Kodaira de la fibre orbifolde $(Y_b,\Delta_b)$ vérifie :
$$\kod{Y_b,\Delta_b}=0.$$
Pour cela, observons que la restriction (encore notée $f$)
$$f:(X_z,D_z)\To (Y_b,\Delta_b)$$
est un morphisme orbifolde entre orbifoldes de même dimension (avec $h(z)=b$). Nous obtenons la conclusion souhaitée :
$$0\leq \kod{Y_b}\leq\kod{Y_b,\Delta_b}\leq \kod{X_z,D_z}=0,$$
la première inégalité étant une conséquence du fait que $Y_b$ est de dimension d'Albanese maximale.
\end{demo}

\section{Analogues anti-hyperboliques}

 Il est conjecturé dans \cite{Ca04} que les propri\'et\'es suivantes (pour $X$ compacte K\"ahler connexe) sont équivalentes:
\begin{enumerate}
\item $X$ est sp\'eciale.
\item La pseudo-m\'etrique de Kobayashi $d_X$ est identiquement nulle.
\item $X$ est $\CC$-connexe (ie: deux quelconques de ses points peuvent \^etre joints par une chaine de courbes enti\`eres trac\'ees dans $X$).
\item Il existe une courbe enti\`ere dense dans $X$.
\item Il existe une courbe enti\`ere Zariski-dense dans $X$.
\item La pseudo-m\'etrique infinit\'esimale $d^*_X$ de Kobayashi est identiquement nulle (cette dernière propriété n'étant pas explicitement considérée dans \emph{loc. cit.}).
\end{enumerate}

Consid\'erons maintenant une vari\'et\'e K\"ahl\'erienne compacte $Y$ ayant m\^eme rev\^etement universel que $X$.

Il est clair que $Y$ poss\`ede les propri\'et\'es (3) et (6) respectivement si $X$ les poss\`ede. Les conjectures pr\'ec\'edentes de \cite{Ca04} impliquent qu'il devrait en \^etre de m\^eme pour les autres propri\'et\'es. Ceci est cependant loin d'\^etre \'evident sans hypoth\`eses additionnelles. Nous faisons quelques remarques \`a ce sujet. 

Supposons par exemple que $\pi_1(X)$ soit presque ab\'elien. Alors $Y$ poss\`ede la propri\'et\'e (1) si c'est le cas pour $X$ d'après le théorème \ref{stab rev universel}. Il serait int\'eressant de savoir si les propri\'et\'es (2),(4) et (5) respectivement pour $X$ l'impliquent alors pour $Y$.

Cette question nous a \'et\'e pos\'ee par Erwan Rousseau (dans le cas (5)) et a motivée l'addition de la présente section.

Nous faisons quelques observations \'el\'ementaires dans le cas des propri\'et\'es (2) et (5), fournissant des indications tr\`es partielles pour leur \'etude, qui montrent la difficult\'e du probl\`eme.

Nous avons tout d'abord un analogue exact du théorème \ref{alb surjective} pour la propri\'et\'e (5). 

\begin{thm}\label{albs} Soit $X$ une vari\'et\'e projective admettant une courbe enti\`ere Zariski dense. Son morphisme d'Albanese est surjectif, \`a fibres connexes et n'a pas de fibre multiple en codimension $1$. 
\end{thm} 
\begin{demo}
C'est une conséquence immédiate de \cite{NWY}, \cite{LW} et \cite{K81}.
\end{demo}

Ces informations ne permettent cependant pas de conclure directement que la courbe enti\`ere Zariski dense $C$ de $X$ relev\'ee \`a son rev\^etement universel redescend en une courbe enti\`ere Zariski-dense sur $Y$ : le choix de $X$ et $Y$ des tores complexes ad\'equats et o\`u $C$ est lin\'eaire montre que ceci est faux en g\'en\'eral. Il semble n\'ecessaire de construire de nouvelles courbes enti\`eres pour pouvoir conclure.

Par contraste, la conclusion de \ref{albs} est ouverte si l'on suppose $d_X=0$ (la surjectivit\'e seule est claire). Le lemme \ref{duniv} permettrait de montrer que la propri\'et\'e (2) se transmet de $X$ \`a $Y$ en faisant les hypoth\`eses additionnelles suivantes sur $X$ :
\begin{itemize}
\item $\pi_1(X)$ est presque ab\'elien,
\item  le morphisme d'Albanese de $X$ est à fibres connexes avec pseudo-m\'etrique de Kobayashi nulle, et sans fibre multiple\footnote{L'assertion la plus difficile \`a \'etablir \'etant, tout comme pour le cas sp\'ecial, l'annulation de la pseudo-m\'etrique de Kobayashi des fibres. Il semble plus int\'er\'essant de tenter de r\'epondre \`a la question g\'en\'erale \ref{d*=d?}.} en codimension $1$.
\end{itemize}

\begin{lem}\label{duniv}
Soit $f:X\to Y$ une fibration entre deux vari\'et\'es projectives complexes lisses et connexes. Soit $\rho:Y'\to Y$ un rev\^etement \'etale (connexe) et $\rho':X':=X\times_Y Y'\to X$ le rev\^etement \'etale d\'eduit par le changement de base $\rho$. De même, notons $F:X'\to Y'$ la fibration déduite de $f$. Si $f$ n'a pas de fibre multiple en codimension $1$, nous avons alors :
\begin{enumerate}[$(i)$]
\item pour tous points $x$ et $y$ dans $X$, $d_X(x,y)\leq \diam(X/Y) + d_Y(f(x),f(y))$ si $\diam(X/Y)$ est le maximum des diam\`etres des fibres de $f$ pour leurs m\'etriques de Kobayashi propres.
\item de même, $d_{X'}(x',y')\leq \diam(X/Y)+d_{Y'}(F(x'),F(y'))$ pour tous points $x'$ et $y'$ dans $X'$.
\item si $\diam(X/Y)=0$ et si $d_{Y'}=0$, alors $d_{X'}$ et $d_X$ sont identiquement nulles.
\end{enumerate}
\end{lem}
\begin{demo}
La démonstration du point $(ii)$ est identique à celle du premier point et le dernier point est une conséquence immédiate des deux premiers : nous démontrons donc le premier point. Soit $V\subset Y$ l'ouvert de Zariski au-dessus duquel les fibres de $f$ sont non-multiples. Sa pseudo-distance de Kobayashi coïncide donc avec la restriction de celle de $Y$ par l'hypoth\`ese de codimension $2$ ou plus (voir par exemple \cite[Th. 3.2.19]{livreKob}. Si $x$ et $y$ sont deux points de $f^{-1}(V)$, nous pouvons joindre $f(x)$ et $f(y)$ par une chaîne de disques holomorphes $\gamma$ contenue dans $V$. En relevant cette chaîne (ce que nous pouvons faire puisque au dessus de $V$ il n'y a pas de fibre multiple), nous en obtenons une, notée $\alpha$, dans $X$ qui joint $x$ à $z$ qui est un point de la fibre $f^{-1}(f(y))$ et nous en déduisons :
\begin{align*}
d_X(x,y)&\le d_X(x,z)+d_X(z,y)\le \longueur(\alpha) +d_{f^{-1}(f(y))}(z,y)\\
&\le \longueur(\gamma)+\diam(X/Y)
\end{align*}
En passant à la limite sur la longueur de $\gamma$, on obtient l'inégalité souhaitée pour $x$ et $y$ dans $f^{-1}(V)$. Par continuité de $d_X$, cette inégalité est donc valable sur $X$ tout entier.
\end{demo}

Les conjectures qui pr\'ec\'edent impliquent une r\'eponse positive \`a la:

\begin{qt}\label{d*=d?}
Soit $X$ K\"ahler compacte telle que $d_X\equiv 0$. A-t-on aussi $d_X^*\equiv 0$?
\end{qt}

\appendix
\section{Résumé orbifolde}
Dans cet appendice, nous pr\'esentons quelques définitions et propriétés des morphismes orbifoldes utilisées ci-dessus. Nous renvoyons \`a  \cite{Ca04,Ca11j} pour une étude plus détaillée de ces notions.

\subsection{La catégorie des orbifoldes : objets et morphismes}
\begin{defi}\label{defi orbi}
Une orbifolde est la donnée d'un couple $(X,\Delta_X)$ où $X$ est une variété complexe et $\Delta_X$ un $\QQ$-diviseur dont les {\bf coefficients} sont compris entre 0 et 1. Les {\bf multiplicités} de $\Delta_X$ sont définies par :
$$\Delta_X=\sum_{D\subset X}(1-\frac{1}{m(D)}).D=\sum_{D\subset X}c(D).D$$
où la somme (localement finie) porte sur tous les diviseurs premiers de $X$. En particulier, $m(D)\ge1$ et $m(D)=1$ pour (localement) presque tout $D$.

Nous dirons que l'orbifolde $(X,\Delta_X)$ est lisse si $X$ l'est et si le support de $\Delta_X$ est à croisements normaux simples.
\end{defi}

Le fait de considérer une telle structure additionnelle est motivée par l'introduction de la base orbifolde d'une fibration (qui tient compte des fibres multiples de cette dernière).
\begin{defi}\label{multiplicité le long d'un div}
Soit $f:X\To Y$ une application holomorphe et $E\subset X$ un diviseur premier de $X$. Nous définissons $m(f,E)$ la multiplicité de $f$ le long de $E$ de la façon suivante :
\begin{enumerate}[(i)]
\item si $E$ est $f$-exceptionnel (ie: si $codim_Y(f(D))\geq 2$) ou si $f(E)=Y$, alors $m(f,E)=1$.
\item sinon $f(E)=D$ est un diviseur premier de $Y$ et nous pouvons écrire
$$f^*(D)=m(f,E)E+F$$
où $F$ est un diviseur effectif de $X$ ne contenant pas $E$ dans son support.
\end{enumerate}
\end{defi}
\noindent Plus généralement :
\begin{defi}\label{base orbi}
Soit $f:X\To Y$ une application holomorphe et supposons $X$ munie d'une structure orbifolde $\Delta_X$. Pour tout diviseur premier $D\subset Y$, nous pouvons écrire :
$$f^*(D)=\sum_{E\subset X\vert f(E)=D} m(f,E)E+R$$
où $R$ est $f$-exceptionnel. Le diviseur
$$\Delta(f,\Delta_X):=\sum_{D\subset Y}(1-\frac{1}{m_f(D)})D$$
est appelé diviseur orbifolde de $f$ (induit par $\Delta_X$) avec :
$$m_f(D):=\mathrm{inf}\{m(f,E).m_X(E)\vert f(E)=D\}.$$
\end{defi}

La catégorie des orbifoldes lisses est munie de la notion de morphismes : ce sont des applications holomorphes entre les variétés sous-jacentes aux orbifoldes et qui doivent respecter les structures additionnelles.
\begin{defi}\label{defi morphisme orbi}
Soit $f:X\To Y$ une application holomorphe et supposons $X$ et $Y$ munies de structures orbifoldes $\Delta_X$ et $\Delta_Y$ respectivement telles que $(X,\Delta_X)$ et $(Y,\Delta_Y)$ soient lisses. Nous dirons que $f$ est un morphisme orbifolde entre $(X,\Delta_X)$ et $(Y,\Delta_Y)$ si, pour tout diviseur premier $D\subset Y$, nous avons :
$$m_i .m_X(E_i)\ge m_Y(D)$$
où nous avons écrit la décomposition irréductible de $f^*(D)$ sous la forme :
$$f^*(D)=\sum_i m_iE_i.$$
Si de plus $f$ est biméromorphe et si $f_*(\Delta_X)=\Delta_Y$, $f$ sera alors appelé un morphisme biméromorphe orbifolde.
\end{defi}
\begin{rem}\label{rem base orbi}
Si $f:X\To Y$ est une application holomorphe et si $\Delta_X$ est une structure orbifolde sur $X$, le morphisme $f:(X,\Delta_X)\To (Y,\Delta(f,\Delta_X))$ n'est en général pas un morphisme orbifolde puisque les diviseurs $f$-exceptionnels ne sont pas pris en compte dans la définition de $\Delta(f,\Delta_X)$.
\end{rem}

\subsection{Différentielles symétriques orbifoldes et fibrations nettes}

On se place dans la situation d'une orbifolde lisse $(X,\Delta_X)$. Les différentielles symétriques adaptées à la structure orbifolde se définissent de la manière suivante :
\begin{defi}\label{defi diff symétrique orbi}
Soit $(x_1,\dots,x_n)$ un système de coordonnées locales dans lequel le diviseur orbifolde a pour équation (symbolique) :
$$\Delta_X=\left(\prod_{i=1}^n x_i^{\big(1-\dfrac{1}{m_i}\big)}=0\right).$$
Le faisceau $\mathbf{S}^m_q(X,\Delta_X)$ est le faisceau de $\mathcal{O}_X$-modules localement libre engendré par les éléments suivants :
$$\bigotimes_{j=1}^m x_j^{\lceil -k_ja_j\rceil}\mathrm{d}x_{J_j}.$$
Dans l'écriture ci-dessus, $\lceil a\rceil$ désigne sa partie entière supérieure et on a de plus :
\begin{enumerate}
\item $a_j=1-\dfrac{1}{m_j}$ est le coefficient de $(x_j=0)$ dans $\Delta_X$.
\item $J_1,\dots,J_m$ est une suite de parties ordonnées à $q$ éléments de l'ensemble des indices $\{1,\dots,n\}$.
\item pour tout $j=1\dots n$, on note $k_j$ le nombre d'occurrences de l'indice $j$ dans la suite $J_1,\dots,J_m$.
\end{enumerate}
\end{defi}
\begin{rem}\label{rem diff symétrique rang max}
Lorsque $q=n=\dim(X)$, nous avons bien entendu
$$\mathbf{S}^m_n(X,\Delta_X)=\lfloor m(K_X+\Delta_X)\rfloor,$$
et dans le cas logarithmique (o\`u tous les $c_j=1$, ou encore: $m_j=+\infty$), nous retrouvons les diff\'erentielles logarithmiques usuelles et leurs puissances sym\'etriques.
\end{rem}
Les différentielles ci-dessus définies sont fonctorielles pour les morphismes orbifoldes :
\begin{prop}\label{formes diff et morphisme}
Soit $f:(X,\Delta_X)\To (Y,\Delta_Y)$ un morphisme orbifolde (entre orbifoldes lisses). L'image réciproque par $f$ des différentielles orbifoldes est bien définie :
$$f^*:\mathbf{S}^m_q(Y,\Delta_Y)\To \mathbf{S}^m_q(X,\Delta_X).$$
\end{prop}

Bien que la base orbifolde d'un morphisme ne suffise pas à rendre ce morphisme un morphisme orbifolde (\emph{cf.} remarque \ref{rem base orbi}), il est cependant possible de prolonger les images réciproques des tenseurs holomorphes dont les pôles sont contrôlés par la base orbifolde, une fois choisi un bon modèle de la fibration initiale. 
\begin{defi}\label{defi fibration nette}
Une fibration $f:(X,\Delta_X)\To Y$ (avec $(X,\Delta_X)$ lisse) est dite nette si il existe $(Z,\Delta_Z)$ une orbifolde lisse et $u:(X,\Delta_X)\To (Z,\Delta_Z)$ un morphisme biméromorphe orbifolde tels que tout diviseur $f$-exceptionnel soit également $u$-exceptionnel.
\end{defi}

Cette notion est plus forte que celle utilis\'ee dans \cite{Ca11j}, mais suffit ici.
\begin{prop}\label{existence modele net}
Soit $g:(Z,\Delta_Z)\To W$ une fibration avec $(Z,\Delta_Z)$ lisse. Il existe alors un diagramme
$$\xymatrix{(X,\Delta_X)\ar[r]^u \ar[d]_f & (Z,\Delta_Z)\ar[d]^g\\
Y\ar[r]^v & W
}$$
dans lequel $u$ et $v$ sont biméromorphes (au sens orbifolde pour $u$) et $f$ est nette.
\end{prop}
Les formes différentielles orbifoldes sont fonctorielles pour les fibrations nettes.
\begin{prop}\label{relevement base orbi nette}
Soit $f:(X,\Delta_X)\To Y$ une fibration nette et considérons $\Delta_Y:=\Delta(f,\Delta_X)$ la base orbifolde de la fibration $f$. Il existe un morphisme de faisceaux bien défini pour $m\ge1$ suffisamment divisible :
$$f^*:\mathcal{O}_Y\left(m(K_Y+\Delta_Y)\right)\To\mathbf{S}^m_q(X,\Delta_X).$$
En particulier, au niveau des sections globales, nous en déduisons l'injection :
$$f^*:H^0\left(Y,\mathcal{O}_Y(m(K_Y+\Delta_Y))\right)\hookrightarrow H^0\left(X,\mathbf{S}^m_q(X,\Delta_X)\right).$$
Si de plus $g:Z\To X$ est une fibration nette telle que $\Delta(g)\ge\Delta_X$, alors on a même :
$$(f\circ g)^*:H^0\left(Y,\mathcal{O}_Y(m(K_Y+\Delta_Y))\right)\hookrightarrow H^0\left(Z,\mathbf{S}^m\Omega^q_Z\right).$$
\end{prop}

\

\subsection{Orbifoldes spéciales}

Les orbifoldes spéciales sont celles qui n'admettent pas de fibration nette sur une base de type général. Plus précisément :
\begin{defi}\label{defi orbi spéciale}
Une orbifolde $(X,\Delta_X)$ est dite spéciale si, pour tout $u:(X',\Delta')\To (X,\Delta_X)$ biméromorphe (au sens de la définition \ref{defi morphisme orbi}) et tout fibration nette $f:(X',\Delta')\To Y$ avec $\dim(Y)>0$, l'inégalité suivante est satisfaite :
$$\kappa(Y,K_Y+\Delta(f,\Delta'))<\dim(Y).$$
\end{defi}

Il est également possible de caractériser les orbifoldes spéciales par leurs faisceaux de formes différentielles. Nous renvoyons à \cite{Ca11j} pour le cas $\Delta\neq 0$ de la caractérisation ci-dessous et nous nous contentons d'énoncer le cas des variétés.
\begin{prop}\label{prop caracterisation speciale}
Une variété $X$ est spéciale si et seulement si pour tout sous faisceau cohérent de rang un $\mathcal{L}\subset \Omega^p_X$ (pour $1\le p\le\dim(X)$), l'inégalité suivante est vérifiée :
$$\kappa(X,\mathcal{L})<p.$$
\end{prop}
%\bibliographystyle{amsalpha}
%\bibliography{bib_speciale}

\begin{thebibliography}{CCE13}

\bibitem[BC13]{BC13}
Caucher Birkar and Jungkai~Alfred Chen, \emph{Varieties fibred over abelian
  varieties with fibres of log general type}, Preprint arXiv:1311.7396, 2013.

\bibitem[C98]{ca98}
Fr{\'e}d{\'e}ric Campana, \emph{Connexit\'e ab\'elienne des vari\'et\'es
  k\"ahl\'eriennes compactes}, Bull. Soc. Math. France \textbf{126} (1998),
  no.~4, 483--506.

\bibitem[C04]{Ca04}
\bysame, \emph{Orbifolds, special varieties and classification theory}, Ann.
  Inst. Fourier (Grenoble) \textbf{54} (2004), no.~3, 499--630.

\bibitem[C11]{Ca11j}
\bysame, \emph{Orbifoldes g\'eom\'etriques sp\'eciales et classification
  bim\'eromorphe des vari\'et\'es k\"ahl\'eriennes compactes}, J. Inst. Math.
  Jussieu \textbf{10} (2011), no.~4, 809--934.

\bibitem[CC14]{CC14}
Fr{\'e}deric Campana and Beno{\^{\i}}t Claudon, \emph{Abelianity conjecture for
  special compact {K}\"ahler 3-folds}, Proc. Edinb. Math. Soc. (2) \textbf{57}
  (2014), no.~1, 55--78.

\bibitem[CCE13]{CCE}
Fr{\'e}d{\'e}ric Campana, Beno{\^i}t Claudon, and Philippe Eyssidieux,
  \emph{Repr{\'e}sentations lin{\'e}aires des groupes k{\"a}hl{\'e}riens :
  Factorisations et conjecture de shafarevich lin{\'e}aire}, Preprint
  arXiv:1302.5016, to appear in Comp.ositio Math., 2013.

\bibitem[CH11]{CH11}
Jungkai~Alfred Chen and Christopher~D. Hacon, \emph{Kodaira dimension of
  irregular varieties}, Invent. Math. \textbf{186} (2011), no.~3, 481--500.

\bibitem[Kaw81]{K81}
Yujiro Kawamata, \emph{Characterization of abelian varieties}, Compositio Math.
  \textbf{43} (1981), no.~2, 253--276.

\bibitem[KO75]{KO75}
Shoshichi Kobayashi and Takushiro Ochiai, \emph{Meromorphic mappings onto
  compact complex spaces of general type}, Invent. Math. \textbf{31} (1975),
  no.~1, 7--16.

\bibitem[Kob98]{livreKob}
Shoshichi Kobayashi, \emph{Hyperbolic complex spaces}, Grundlehren der
  Mathematischen Wissenschaften [Fundamental Principles of Mathematical
  Sciences], vol. 318, Springer-Verlag, Berlin, 1998.
  
\bibitem[LW12]{LW}
Steven Lu and J\" org Winkelmann, \emph{Quasiprojective varieties admitting Zariski dense entire holomorphic curves}, Forum Math. \textbf{24} (2012), no.~2, 399--418.


\bibitem[NWY07]{NWY}
Juniro Noguchi, J\" org Winkelmann and Katsutochi Yamanoi, \emph{Degeneracy of holomorphic curves into algebraic varieties},  J. Math. Pures Appl. (9) \textbf{88} (2007), no.~3, 293--306.

\bibitem[Sak74]{Sak}
Fumio Sakai, \emph{Degeneracy of holomorphic maps with ramification}, Invent.
  Math. \textbf{26} (1974), 213--229.

\bibitem[Siu98]{S98}
Yum-Tong Siu, \emph{Invariance of plurigenera}, Invent. Math. \textbf{134}
  (1998), no.~3, 661--673.

\bibitem[Siu02]{S02}
\bysame, \emph{Extension of twisted pluricanonical sections with
  plurisubharmonic weight and invariance of semipositively twisted plurigenera
  for manifolds not necessarily of general type}, Complex geometry
  ({G}\"ottingen, 2000), Springer, Berlin, 2002, pp.~223--277.

\bibitem[Tsu96]{T96}
Hajime Tsuji, \emph{On the universal covering of projective manifolds of
  general type}, Kodai Math. J. \textbf{19} (1996), no.~1, 137--143.

\bibitem[Uen75]{U75}
Kenji Ueno, \emph{Classification theory of algebraic varieties and compact
  complex spaces}, Lecture Notes in Mathematics, Vol. 439, Springer-Verlag,
  1975, Notes written in collaboration with P. Cherenack.

\end{thebibliography}

\providecommand{\bysame}{\leavevmode\hbox to3em{\hrulefill}\thinspace}
\providecommand{\MR}{\relax\ifhmode\unskip\space\fi MR }
% \MRhref is called by the amsart/book/proc definition of \MR.
\providecommand{\MRhref}[2]{%
  \href{http://www.ams.org/mathscinet-getitem?mr=#1}{#2}
}
\providecommand{\href}[2]{#2}

\end{document}